\documentclass[11 pt,a4paper,twoside,reqno]{amsart} 
\usepackage{amsfonts,amssymb,amscd,amsmath,enumerate,verbatim,calc} 
 
\usepackage{amsmath}
\usepackage{amsfonts}
\usepackage{amssymb}
\usepackage{setspace}
\usepackage{amsthm,mathrsfs,textcomp}
\usepackage[utf8]{inputenc}	
\usepackage[T1]{fontenc}	
\usepackage{mathtools}
\usepackage[super]{nth}

\theoremstyle{plain}

\theoremstyle{definition}

\numberwithin{equation}{section}

\newcommand {\Mathstu}{\hspace*{-0.6 cm}The Mathematics Student\hfill \hspace{2 in}\hfill ISSN: 0025-5742\\ Vol. 91, Nos. 3-4, July-December (2022), xx--yy}
\begin{document}
	\Mathstu
	\vspace{0.7 cm}
	
	\title
	[Solution of single parameter Bring quintic equation]
	{\scshape \large \bf {Solution of single parameter Bring quintic equation}}
	\vspace*{-0.4 cm}
	\author [Raghavendra G. Kulkarni] {Raghavendra G. Kulkarni}
	\thanks{\raggedright{\hspace*{-.21 in}
\hspace{0.24cm}			{2010 Mathematics Subject Classification}: 12E12, 40A05    }  \\ [0.4 cm]
\hspace{5.5cm}		\raggedleft \textit{\large \copyright\,\,{Indian Mathematical Society, 2021}}}
	
	\maketitle
	\vspace*{-0.9 cm}

\begin{abstract}
In this paper, we propose a new method to obtain a solution to a single-parameter Bring quintic equation of the form, $x^{5}+x=a$, where $a$ is real. The method transforms the given quintic equation to an infinite but convergent series expression in $(x/a)$, which is further transformed to a quartic equation in a novel fashion. The coefficients of the quartic equation so obtained are some kind of infinite series expressions in $a^{-4}$, which are termed as \textit{ultraradicals}. The quartic equation is then solved and its one real solution is picked; further using this, the real solution of quintic equation, $x^{5}+x=a,$ is extracted. The ultraradicals used in this method converge for $|a| > 1$; hence the method can be used when $|a| > 1$.
\end{abstract}

\maketitle

\section{Introduction}
Many real world applications encounter quintic equations. For example, the three-body motion of celestial objects governed by Keppler laws requires solution to a quintic equation to determine the position of the objects \cite{NI85}. In hydraulic engineering, the height of water flowing in an open rectangular channel is obtained by solving a quintic equation \cite{KM11}. In structural mechanics, the study of nonlinear vibration of beams involves quintic equations \cite{HM13}. Since the general quintic equations do not have closed form solutions (similar to that exist for cubic and quartic equations), one has to resort to either numerical methods or use special functions (called \textit{ultraradicals}) to obtain the solutions. A brief discussion on the historical account of efforts put in by several mathematicians for solving quintic equations is given below.
\\

The cubic and the quartic equations were solved by Cardan and Ferrari respectively in the sixteenth century; and the intense struggle by several mathematicians to solve the quintic equation, using the techniques similar to that adopted for cubic and quartic equations met with no success. In 1683, the German mathematician Ehrenfried Walther von Tschirnhaus introduced a polynomial transformation bearing his name, which could eliminate intermediate terms in an $n$-th degree polynomial equation. Using the transformation, Tschirnhaus showed that the cubic equation can be solved by removing the two intermediate terms. Further, the quartic equation was solved by removing the two odd-power terms \cite{JB06}. This result encouraged the mathematicians to explore the ways for reducing a polynomial equation to simpler forms, hoping that it might lead to the solutions of the quintic and the higher degree equations.
\\

In 1786, the Swedish mathematician E. S. Bring showed that one can transform the general quintic equation to the form $x^{5}+Ax+B=0$. It appears that the papers containing the works of Bring got lost in the archives of University of Lund. Probably unaware of Bring's work, the French mathematician G. B. Jerrard (1852) also reduced the general quintic equation to the form $x^{5}+Ax+B=0$ \cite{JB06,VA03}; for detailed description of getting this form using Tschirnhaus transformation, see \cite{LE26,GB59}.
\\

However, all these efforts of reducing the polynomial equations to simpler forms could not achieve the desired goal of solving the quintic and the higher degree equations. The issue of unsolvability of these equations was eventually resolved by Abel (1824) and Galois (1832), who conclusively proved that the general polynomial equations of degree five and above cannot be solved by the methods similar to that were adopted for solving the cubic and the quartic equations \cite{JB06}. This means that these equations cannot be solved by performing finite number of elementary arithmetic operations (addition, subtraction, multiplication, and division) along with the radical ($n$-th root) operation, which were enough to solve the quadratic, the cubic, or the quartic equations. In other words, quintic and the higher-degree equations cannot be solved using radicals alone.
\\ 

The simplest form of quintic equation that is unsolvable in radicals is Bring-Jerrard quintic equation, $x^{5}+Ax+B=0$. This doesn't mean that there are no solutions to the quintic equations; the solutions of quintics can still be obtained (without violating the Galois theory) by the use of more complicated functions than the radicals. Notice that the radical ($n$-th root, $ \; n=2,3,4$) operation can be viewed as a transcendental function involving a logarithmic or trigonometric function in a parameter \cite{BK96}. The complicated functions that are useful for solving the quintic and the higher-degree equations are now termed as \textit{ultraradicals} \cite{BK96}. Elliptic modular functions are one form of ultraradicals used by Hermite (1858) for solving the single-parameter Bring quintic equation, $x^{5}-x-a=0$ \cite{RB24,MK15}.
\\

Also we notice from the literature that the Bring quintic equation of the form, $x^{5}+x=a$ (with $a$ real), is solved using an ultraradical called Bring radical, $BR(a)$, which is an infinite series function in $a$, as shown below \cite{GB59,WIKI}.
\[ x=BR(a)=a-a^{5}+5a^{9}-35a^{13}+285a^{17}-2530a^{21}+ \ldots \]
The above series is convergent for values of $|a|< 4/( \; 5 \cdot \sqrt [4] {5} \; )$, see \cite{WIKI}.
\\  

In this paper, we present a method for extracting a real solution of the quintic equation, $x^{5}+x=a$ (with $a$ real), by transforming it to a quartic equation, whose coefficients are ultraradicals. The ultraradicals in our case happen to be some kind of infinite series functions in $a^{-4}$, as one can see later in this paper. These series functions are convergent for $|a| > 1$; so the proposed method can be used when $|a| > 1$.

\section{The proposed method}
Consider the following single-parameter Bring quintic equation,
\begin{equation}
x^{5}+x=a,
\end{equation}
where $a$ is real. Since any polynomial of odd-degree in $x$ (with real coefficients) crosses the $x$-axis at least once, the quintic equation (2.1) has at least one real solution. We plan to extract this solution algebraically. For this purpose, when (2.1) is rearranged as, $(x/a)=1/(1+x^{4}),$ it is easier to notice that the real root ($x$) of quintic equation (2.1) satisfies the inequality, $0 < (x/a) < 1$, since $x^{4}$ is always a positive number. Equipped with this information, (2.1) is further rearranged as,
\begin{equation}
(x/a)^{5}=a^{-4}[1-(x/a)].
\end{equation}
Taking the fifth root of (2.2) yields,
\begin{equation}
(x/a)=a^{-4/5}[1-(x/a)]^{1/5}.
\end{equation}
Expanding the right-hand-side (RHS) of (2.3) using binomial theorem results in,
\begin{equation}
(x/a)=a^{-4/5} \left[1+\sum_{k=1}^{\infty}\frac{n(n-1)(n-2) \; ...\; k \; \mbox{terms}}{k!}\left(\frac{-x}{a}\right)^{k}\right],
\end{equation}
where $n=1/5$. Expressing the infinite series in (2.4) explicitly and rearranging the terms yields the following expression,
\begin{eqnarray}
a^{4/5}(x/a)=1-c_{1}(x/a)-c_{2}(x/a)^{2}-c_{3}(x/a)^{3}-c_{4}(x/a)^{4}-c_{5}(x/a)^{5} \nonumber\\
-c_{6}(x/a)^{6}-c_{7}(x/a)^{7}-c_{8}(x/a)^{8}-c_{9}(x/a)^{9}-c_{10}(x/a)^{10}- \; ...,
\end{eqnarray}
where '$c_{k}$'s are evaluated using the recurrence formula, 
\[ c_{k+1}=[(5k-1)/5(k+1)]c_{k}, \]
and are listed in Table \ref{table:tab1}.
\begin{table}[h!]
\caption{Values of $c_{k}$ for $k=1$ to $36$}
\centering
\begin{tabular}{r l c c c c}
\hline\hline 
$k$ & $c_{k}$ & $k$ & $c_{k}$ & $k$ & $c_{k}$ \\ [0.5ex]
\hline
1 & 0.2 & 13 & 0.007986216075264 & 25 & 0.003627166639633  \\
2 & 0.08 & 14 & 0.007301683268813 & 26 & 0.003459758948573  \\
3 & 0.048 & 15 & 0.006717548607308 & 27 & 0.003305991884192 \\
4 & 0.0336 & 16 & 0.00621373246176 & 28 & 0.003164306517727 \\
5 & 0.025536 & 17 & 0.005775116052694 & 29 & 0.003033369696303 \\
6 & 0.0204288 & 18 & 0.005390108315848 & 30 & 0.002912034908451 \\
7 & 0.01692672 & 19 & 0.005049680422216 & 31 & 0.002799310976511 \\
8 & 0.014387712 & 20 & 0.004746699596883 & 32 & 0.002694336814892 \\
9 & 0.0124693504 & 21 & 0.004475459619918 & 33 & 0.002596360930714 \\
10 & 0.010973028352 & 22 & 0.00423134364065 & 34 & 0.002504724662571 \\
11 & 0.0097759707136 & 23 & 0.004010577885485 & 35 & 0.002418848388426 \\
12 & 0.00879837364224 & 24 & 0.003810048991211 & 36 & 0.002338220108812 \\ [0.5ex]
\hline\hline
\end{tabular}
\label{table:tab1}
\end{table}
Since $0 < (x/a) < 1$ and $0 < c_{k+1} < c_{k} < 1$, it implies that the terms in the infinite series in (2.5) satisfy the inequality, $c_{k+1}(x/a)^{k+1} < c_{k}(x/a)^{k}$; hence we conclude that the series is a convergent infinite series. Rearranging (2.5) yields,
\begin{eqnarray}
1-(c_{1}+a^{4/5})(x/a)-c_{2}(x/a)^{2}-c_{3}(x/a)^{3}-c_{4}(x/a)^{4}-c_{5}(x/a)^{5} \nonumber\\
-c_{6}(x/a)^{6}-c_{7}(x/a)^{7}-c_{8}(x/a)^{8}-c_{9}(x/a)^{9}-...=0,
\end{eqnarray}
Now we use (2.2) repeatedly to eliminate $(x/a)^{5}$ and higher power terms from (2.6), which results in a quartic equation in $(x/a)$ as below.
\begin{equation}
K_{4}(x/a)^{4}+K_{3}(x/a)^{3}+K_{2}(x/a)^{2}+K_{1}(x/a)+K_{0}=0,
\end{equation}
where $K_{0}$, $K_{1}$, $K_{2}$, $K_{3}$, and $K_{4}$ are given by:
\begin{eqnarray}
K_{0}=1-(c_{5}/a^{4})+[(c_{9}-c_{10})/a^{8}]-[(c_{13}-2c_{14}+c_{15})/a^{12}]+ \; ... \nonumber\\
K_{1}=-a^{4/5}-c_{1}+[(c_{5}-c_{6})/a^{4}]-[(c_{9}-2c_{10}+c_{11})/a^{8}] \nonumber\\
+[(c_{13}-3c_{14}+3c_{15}-c_{16})/a^{12}]- \; ... \nonumber\\
K_{2}=-c_{2}+[(c_{6}-c_{7})/a^{4}]-[(c_{10}-2c_{11}+c_{12})/a^{8}] \nonumber\\
+[(c_{14}-3c_{15}+3c_{16}-c_{17})/a^{12}]- \; ... \nonumber\\
K_{3}=-c_{3}+[(c_{7}-c_{8})/a^{4}]-[(c_{11}-2c_{12}+c_{13})/a^{8}] \nonumber\\
+[(c_{15}-3c_{16}+3c_{17}-c_{18})/a^{12}]- \; ... \nonumber\\
K_{4}=-c_{4}+[(c_{8}-c_{9})/a^{4}]-[(c_{12}-2c_{13}+c_{14})/a^{8}] \nonumber\\
+[(c_{16}-3c_{17}+3c_{18}-c_{19}/a^{12}]+ \; ... 
\end{eqnarray}
Notice that the coefficients, $K_{0}$, $K_{1}$, $K_{2}$, $K_{3}$, and $K_{4}$, given in (2.8) are infinite series expressions in $a^{-4}$, and are termed as ultraradicals in $a$, similar to that defined in \cite{BK96}. The convergence of these functions will be discussed in the next section. Normalizing the quartic equation (2.7) by dividing it throughout by $K_{4}$ yields,
\begin{eqnarray}
(x/a)^{4}+(K_{3}/K_{4})(x/a)^{3}+(K_{2}/K_{4})(x/a)^{2} \nonumber\\
+(K_{1}/K_{4})(x/a)+(K_{0}/K_{4})=0.
\end{eqnarray}
Solving (2.9) by traditional methods \cite{LE22,RG13} yields four solutions of $(x/a)$. Since we already know that a real solution of the given quintic equation (2.1) satisfies the inequality, $0 < (x/a) < 1$, we now pick a solution of (2.9), which is real, positive, and less than unity, say $x_{1}/a$; from this the desired solution of (2.1) is obtained as $x_{1}$.

\section{Salient features of ultraradicals}
In order to study the salient features of the ultraradicals defined above, we first express them in general forms as below.
\begin{eqnarray}
K_{0}=1+\sum_{m=1}^{\infty}\left[\frac{(-1)^{m}}{a^{4m}}\sum_{n=0}^{m-1}(-1)^{n}C_{n}^{m-1}c_{4m+n+1} \right], \nonumber\\
K_{1}=-a^{4/5}+\sum_{m=1}^{\infty}\left[\frac{(-1)^{m}}{a^{4(m-1)}}\sum_{n=0}^{m-1}(-1)^{n}C_{n}^{m-1}c_{4m+n-3}\right],  \nonumber\\
K_{2}=\sum_{m=1}^{\infty}\left[\frac{(-1)^{m}}{a^{4(m-1)}}\sum_{n=0}^{m-1}(-1)^{n}C_{n}^{m-1}c_{4m+n-2}\right], \nonumber\\
K_{3}=\sum_{m=1}^{\infty}\left[\frac{(-1)^{m}}{a^{4(m-1)}}\sum_{n=0}^{m-1}(-1)^{n}C_{n}^{m-1}c_{4m+n-1}\right], \nonumber\\
K_{4}=\sum_{m=1}^{\infty}\left[\frac{(-1)^{m}}{a^{4(m-1)}}\sum_{n=0}^{m-1}(-1)^{n}C_{n}^{m-1}c_{4m+n}\right].
\end{eqnarray}
Further, making use of numerical values given in Table 1, the expressions in (3.1) are written as:
\begin{eqnarray}
K_{0}=1-(0.025536)a^{-4}+(0.001496322...)a^{-8} \nonumber\\
-(0.000100398...)a^{-12}+(0.000007132...)a^{-16}- \ldots, \nonumber\\
K_{1}=-a^{4/5}-0.2\left[1-(0.025536)a^{-4}+(0.001496322...)a^{-8} \right. \nonumber\\
\left. -(0.000100398...)a^{-12}+(0.000007132...)a^{-16}- \ldots \right], \nonumber\\
K_{2}=0.92-\left[1-(0.00350208)a^{-4}+(0.0002194606...)a^{-8} \right. \nonumber\\
\left. -(0.0000151188...)a^{-12}+(0.0000010894...)a^{-16}- \ldots \right], \nonumber\\
K_{3}=0.952-\left[1-(0.002539008)a^{-4}+(0.0001654395...)a^{-8} \right.  \nonumber\\
\left. -(0.0000115911...)a^{-12}+(0.0000008431...)a^{-16}- \ldots \right], \nonumber\\
K_{4}=0.9664-\left[1-(0.0019183616...)a^{-4}+(0.0001276248...)a^{-8} \right. \nonumber\\
\left. -(0.0000090288...)a^{-12}+(0.0000006604...)a^{-16}- \ldots \right].
\end{eqnarray}

It is essential to know for which values of $a$ the series expressions in (3.2) are convergent, so that the proposed method becomes relevant in that range. Observing the expressions given in (3.2), we notice that $K_{1}=-a^{4/5}-0.2K_{0}$. Further, a comparison of infinite series in $K_{2}$ with that of $K_{0}$ reveals that the magnitudes of coefficients of $a^{-4m}$ in $K_{2}$ are smaller than that in $K_{0}$ ($m=1,2,3, \ldots $); hence, if $K_{0}$ is convergent for some value of $a$, $K_{2}$ also will be convergent for that $a$.
\\ 

Similarly the magnitudes of coefficients of $a^{-4m}$ in $K_{3}$ and $K_{4}$ are also smaller than that in $K_{0}$. Hence $K_{3}$ and $K_{4}$ will also be convergent for that value of $a$, for which $K_{0}$ is convergent. Therefore, it is sufficient if we determine the convergence of $K_{0}$, and this can be done more easily with experimental techniques than the theoretical ones. For this purpose, we determine first few terms ($T_{m}$) in $K_{0}$. The first forty-one terms in $K_{0}$ are determined using the expression in (3.1) and  listed in Table \ref{table:tab2}, for $|a|=1$.  The partial sums of the (first) eleven terms ($S_{11}$), twenty-one terms ($S_{21}$), thirty-one terms ($S_{31}$), and the  forty-one terms ($S_{41}$) for various values of $|a|$ are determined and listed in Table \ref{table:tab3}.
\\

\begin{table}[h!]
\caption{First forty-one terms in $K_{0}$ for $|a|=1$.}
\centering
\begin{tabular}{c c c c c c}
\hline\hline 
$m$ & $T_{m}$ & $m$ & $T_{m}$ & $m$ & $T_{m}$ \\ [0.5ex]
\hline
0 & 1 & 14 & -4.18719E-16 & 28 & -5.91874E-13  \\
1 & -0.025536 & 15 & 2.03407E-15 & 29 & 3.10788E-12  \\
2 & 0.001496322 & 16 & -8.04998E-15 & 30 & 4.96923E-12 \\
3 & -0.000100398 & 17 & -3.5667E-15 & 31 & -1.00878E-11 \\
4 & 7.13277E-06 & 18 & 3.38033E-15 & 32 & 3.2447E-12 \\
5 & -5.23045E-07 & 19 & -3.77267E-14 & 33 & -1.26564E-10 \\
6 & 3.91345E-08 & 20 & 3.64555E-14 & 34 & -1.60997E-11 \\
7 & -2.96913E-09 & 21 & 9.31319E-14 & 35 & 1.56419E-10 \\
8 & 2.27581E-10 & 22 & 1.88169E-14 & 36 & -1.14683E-10 \\
9 & -1.75807E-11 & 23 & 4.60115E-14 & 37 & 3.06899E-10 \\
10 & 1.36658E-12 & 24 & 8.33396E-14 & 38 & 9.61151E-10 \\
11 & -1.06755E-13 & 25 & -6.6296E-14 & 39 & 5.34989E-09 \\
12 & 7.92964E-15 & 26 & 6.9709E-13 & 40 & 1.63177E-09 \\    
13 & -1.76335E-15 & 27 & 4.35009E-14 &   &   \\ [0.5ex]
\hline\hline
\end{tabular}
\label{table:tab2}
\end{table}
\begin{table}[h!]
\caption{Partial sums of terms in $K_{0}$ for various values of $|a|$.}
\centering
\begin{tabular}{c c c c c}
\hline\hline 
$|a|$ & $S_{11}$ & $S_{21}$ & $S_{31}$ & $S_{41}$ \\ [0.5ex]
\hline
1.5 & 0.995013473 & 0.995013473 & 0.995013473 & 0.995013473 \\
1.2 & 0.988022294 & 0.988022294 & 0.988022294 & 0.988022294  \\
1 & 0.97586657 & 0.97586657 & 0.97586657 & 0.975866578  \\
0.9 & 0.964234296 & 0.964234296 & 0.964236432 & 1.082230751 \\
0.85 & 0.955950756 & 0.955950764 & 0.957864374 & 928.7584625  \\ [0.5ex]
\hline\hline
\end{tabular}
\label{table:tab3}
\end{table}
From the Table \ref{table:tab3}, we note that for $|a| \geq 1$, all the partial sums ($S_{11}$, $S_{21}$, $S_{31}$, and $S_{41}$) are more or less the same, which means adding more terms to the sum has insignificant effect on the value of the sum, implying the series is convergent. However, for $|a| < 1$, the partial sums progressively become larger ($S_{41} > S_{31} > S_{21} > S_{11}$); so, the sum of infinite number of terms tends to infinity, which means the series is divergent.

\section{Numerical example}
Consider the quintic equation,
\[ x^{5}+x=9.09375, \]
for which one real solution has to be determined using the method proposed here. For this purpose, we obtain the infinite series expression from (2.5), and a polynomial equation from (2.6). The ultraradicals,  $K_{0}$, $K_{1}$, $K_{2}$, $K_{3}$, and $K_{4}$ are determined from (2.9) as:
\begin{eqnarray*}
K_{0}=0.999996266\ldots, \; \; K_{1}=-6.047824804\ldots, \; \; K_{2}=-0.079999488\ldots, \\
K_{3}=-0.047999629\ldots, \; \; K_{4}=-0.033599719\ldots,
\end{eqnarray*}
and the coefficients of quartic equation (3.1) are obtained as follows:
\begin{eqnarray*}
\frac{K_{0}}{K_{4}}=-29.7620421\ldots, \; \; \frac{K_{1}}{K_{4}}=179.9962886\ldots, \\
\frac{K_{2}}{K_{4}}=2.380957018\ldots, \; \; \frac{K_{3}}{K_{4}}=1.428572306\ldots .
\end{eqnarray*}
Using these coefficients, the quartic equation (3.1) is solved using the methods available in literature (see \cite{LE22,RG13}); and the real solution---that is positive and less than unity---is picked as: $(x_{1}/a)=0.1649484536\ldots , $ and the real solution of the given quintic equation, $x^{5}+x=9.09375$, is then obtained as: $x_{1}=1.50000000000009\ldots,$ while the exact solution is $1.5$; the difference between the exact value and the determined value is due to the limited number of terms used in the infinite series in (2.5). In this case we have used 15 terms.
\\

Let us use Newton method for finding a real root of the given polynomial, $f(x)=x^{5}+x-9.09375$, numerically. The Newton method provides an expression as given below for the iterative procedure to obtain a better approximation ($x_{n+1}$) to the root than the previous approximation ($x_{n}$).
\[ x_{n+1}=x_{n}-[f(x_{n})/f'(x_{n})] \]
The initial guess for the root is denoted as $x_{0}$, which has to be a proper guess. In the above example, since we know that the value of the real root is greater than one, let $x_{0}=1$. Using the iterative procedure mentioned above, an approximate real root is obtained in seven iterations, which has an accuracy better than that obtained from the proposed method.
\\

\begin{center}\textsc{Concluding comments} \end{center}

This paper has presented a new method to obtain a real solution of a single-parameter Bring quintic equation, $x^{5}+x=a$, where $a$ is real. The method transforms the quintic equation to an infinite but convergent series, which is further transformed into a quartic equation in a novel fashion, whose coefficients are some kind of infinite series functions (ultraradicals) in $a^{-4}$. Thus the task of solving the quintic equation reduces now to solving a quartic equation. The proposed method can be used for values of $|a| > 1$, as the series functions in the ultraradicals converge for $|a| > 1$.
\\

\noindent \textbf{Acknowledgement:}
The author thanks the administration of PES University for supporting this work. Also, the author acknowledges the valuable comments from the reviewer, which improved the manuscript.

\bigskip

\bigskip

\noindent \textsc{Raghavendra G. Kulkarni\\ Department of Electronics \& Communication Engineering\\ PES University\\ 100 Feet Ring Road, BSK III Stage\\ Bengaluru - 560085, INDIA.}\\
E-mail: raghavendrakulkarni@pes.edu

\end{document}